# COAGULATION–FRAGMENTATION DUALITY, POISSON–DIRICHLET DISTRIBUTIONS AND RANDOM RECURSIVE TREES

By Rui Dong,[1] Christina Goldschmidt[2] and James B. Martin[3]

*University of California, Berkeley, Cambridge University and Oxford University*

In this paper we give a new example of duality between fragmentation and coagulation operators. Consider the space of partitions of mass (i.e., decreasing sequences of nonnegative real numbers whose sum is 1) and the two-parameter family of Poisson–Dirichlet distributions $PD(\alpha, \theta)$ that take values in this space. We introduce families of random fragmentation and coagulation operators $Frag_\alpha$ and $Coag_{\alpha,\theta}$, respectively, with the following property: if the input to $Frag_\alpha$ has $PD(\alpha, \theta)$ distribution, then the output has $PD(\alpha, \theta+1)$ distribution, while the reverse is true for $Coag_{\alpha,\theta}$. This result may be proved using a subordinator representation and it provides a companion set of relations to those of Pitman between $PD(\alpha, \theta)$ and $PD(\alpha\beta, \theta)$. Repeated application of the $Frag_\alpha$ operators gives rise to a family of fragmentation chains. We show that these Markov chains can be encoded naturally by certain random recursive trees, and use this representation to give an alternative and more concrete proof of the coagulation–fragmentation duality.

**1. Introduction.** The subject of this paper is a duality relations for a fragmentation operator and a coagulation operator when applied to certain Poisson–Dirichlet distributions. The idea of duality by time reversal for fragmentation and coagulation is very natural: the opposite of splitting blocks apart is coalescing them. However, demonstrating duality for coagulation and fragmentation processes with desirable properties seems to be a difficult problem and there is no general theory. There are, however, several

Received September 2005; revised May 2006.
[1]Supported in part by NSF Grant DMS-04-05779.
[2]Supported by Pembroke College, Cambridge.
[3]Supported by the CNRS.

beautiful examples where some form of duality does hold; for instance, the additive coalescent of Aldous and Pitman [4] and the Bolthausen–Sznitman [10] coalescent, whose duality properties were discovered by Pitman [18] (see also the discussion in Chapter 5 of [19]).

We work on the space of partitions of mass (i.e., decreasing sequences of nonnegative real numbers whose sum is 1). Fix $0 \leq \alpha < 1$ and $\theta > -\alpha$. Our fragmentation operator takes a size-biased pick from the sequence and splits the chosen block with a $\mathrm{PD}(\alpha, 1 - \alpha)$ random variable. Our coagulation operator generates a $\mathrm{Beta}((1 - \alpha)/\alpha, (\theta + \alpha)/\alpha)$ random variable and coalesces that proportion of the blocks. If the input to the fragmentation operator has $\mathrm{PD}(\alpha, \theta)$ distribution, then its output has $\mathrm{PD}(\alpha, \theta + 1)$ distribution. Moreover, an application of the coagulation operator allows us to go back the other way. This extends a result of Bertoin and Goldschmidt [7], which covered the $\alpha = 0$ case. It also provides a companion set of relations to those of [18] for the Poisson–Dirichlet distributions.

Building a Markov process using the fragmentation operator gives a self-similar fragmentation process of index 1, dislocation measure $\mathrm{PD}(\alpha, 1 - \alpha)$ and erosion coefficient 0, in the terminology of Bertoin [6]. We show that this fragmentation process is naturally embedded in certain random recursive trees (i.e., rooted labeled trees whose labels increase along paths away from the root). These $(\alpha, \theta)$-*recursive trees* can be viewed as nested systems of exchangeable partitions, and their construction is an elaboration of the "Chinese restaurant process" of Dubins and Pitman (see [19]). They can also be regarded as examples of graphs constructed by *preferential attachment* [9, 15, 21], and our results complement those which have previously been obtained in special cases, for example, $\alpha = 1/2, \theta = 0$ [23] or $\alpha = 1/2, \theta = 1/2$ [14].

In Section 2 we collect various definitions and results concerning the family of Poisson–Dirichlet distributions. The fragmentation and coagulation operators are defined precisely in Section 3 and the duality relationship between them is proved using a subordinator representation. The extension to fragmentation and coagulation processes is described in Section 4. In Section 5 we introduce the random recursive tree model and describe how it encodes the fragmentation process. In Section 6 we use this representation to give an alternative and more concrete proof of the duality between the fragmentation and coagulation processes. Finally, in Section 7, we comment on relationships between the recursive tree model and previous representations of the Chinese restaurant process in terms of continuous-time branching processes.

**2. Poisson–Dirichlet distributions.** We will be concerned with properties of the two-parameter Poisson–Dirichlet distribution, introduced in its full generality in [20]. We will first define the $\mathrm{PD}(\alpha, \theta)$ distribution and then mention some of its properties which we will use in the sequel.

DEFINITION 2.1 (*Stick-breaking scheme*). For $0 \leq \alpha < 1$ and $\theta > -\alpha$, let $B_1, B_2, \ldots$ be independent random variables such that $B_n \sim \text{Beta}(1 - \alpha, \theta + n\alpha)$ for all $n \geq 1$. Let

$$\tilde{X}_1 = B_1 \quad \text{and} \quad \tilde{X}_n = (1 - B_1) \cdots (1 - B_{n-1}) B_n \qquad \text{for } n \geq 2.$$

Let $X_1 \geq X_2 \geq \cdots$ be the ranked values of the $\tilde{X}_n$. Then define the $\text{PD}(\alpha, \theta)$ distribution to be the law of the vector $(X_1, X_2, \ldots)$.

The sequence $(\tilde{X}_1, \tilde{X}_2, \ldots)$ is a size-biased permutation of $(X_1, X_2, \ldots)$ and is said to have the Griffiths–Engen–McCloskey $\text{GEM}(\alpha, \theta)$ distribution. In particular, $\tilde{X}_1$ is a size-biased pick from $(X_1, X_2, \ldots)$. The next proposition is a direct consequence of Definition 2.1.

PROPOSITION 2.2 ([20]). *Let $(Y_n, n \geq 1) \sim \text{PD}(\alpha, \theta + \alpha)$ and let $B$ be an independent $\text{Beta}(1 - \alpha, \theta + \alpha)$ random variable. Let the sequence $(X_m, m \geq 1)$ be defined by inserting $B$ into the sequence $((1 - B)Y_n, n \geq 1)$ and reranking. Then $(X_m, m \geq 1)$ has $\text{PD}(\alpha, \theta)$ distribution and $B$ is a size-biased pick from $(X_m, m \geq 1)$.*

There are many representations of the $\text{PD}(\alpha, \theta)$ distribution. For $\alpha = 0$ and $\theta > 0$, we have Kingman's subordinator representation:

PROPOSITION 2.3 ([13]). *Let $\gamma$ be a standard gamma subordinator on the time interval $[0, \theta]$, with ranked jumps $\xi_1 > \xi_2 > \cdots > 0$. Then*

$$\frac{1}{\gamma(\theta)}(\xi_1, \xi_2, \ldots) \sim \text{PD}(0, \theta)$$

*independently of $\gamma(\theta)$, which has a $\text{Gamma}(\theta, 1)$ distribution.*

A related subordinator representation holds for $0 < \alpha < 1$ and $\theta > 0$:

PROPOSITION 2.4 ([20]). *Fix $0 < \alpha < 1$. Let $(\tau(s), s \geq 0)$ be a subordinator with Lévy measure $\alpha x^{-\alpha - 1} e^{-x} dx$. Let $S$ be an independent $\text{Gamma}(\theta/\alpha, \Gamma(1 - \alpha))$ random variable and let $\xi_1 > \xi_2 > \cdots > 0$ be the ranked jumps of the subordinator in the time interval $[0, S]$. Then*

$$\frac{1}{\tau(S)}(\xi_1, \xi_2, \ldots) \sim \text{PD}(\alpha, \theta)$$

*independently of $\tau(S)$, which has a $\text{Gamma}(\theta, 1)$ distribution.*

**3. Dual fragmentation and coagulation operators.** Let $0 \leq \alpha < 1$, $\theta > -\alpha$ and

$$\Delta_\infty^\downarrow = \left\{ x = (x_1, x_2, \ldots) : x_1 \geq x_2 \geq \cdots > 0, \quad \sum_{i=1}^\infty x_i = 1 \right\}.$$

Let $\mathrm{Frag}_\alpha : \Delta_\infty^\downarrow \to \Delta_\infty^\downarrow$ be a random operator which takes a size-biased pick from its input, splits it using an independent $\mathrm{PD}(\alpha, 1-\alpha)$ random variable and then puts the resulting vector in decreasing order. More precisely, fix $x \in \Delta_\infty^\downarrow$. Let $I$ be an index chosen according to the distribution

$$\mathbb{P}(I = i) = x_i, \qquad i \geq 1,$$

and let $\eta = (\eta_1, \eta_2, \ldots) \sim \mathrm{PD}(\alpha, 1-\alpha)$ independently of $I$. Then

$$\mathrm{Frag}_\alpha(x) = (x_1, x_2, \ldots, x_{I-1}, x_I \eta_1, x_I \eta_2, \ldots, x_{I+1}, x_{I+2}, \ldots)^\downarrow,$$

where here, as throughout this paper, the arrow used as a superscript on a sequence means that the sequence is to be put into decreasing order. Let $\mathrm{Coag}_{\alpha,\theta} : \Delta_\infty^\downarrow \to \Delta_\infty^\downarrow$ be another random operator which picks a $\mathrm{Beta}((1-\alpha)/\alpha, (\theta+\alpha)/\alpha)$ proportion of the blocks if $\alpha > 0$, or a deterministic proportion $1/(\theta+1)$ if $\alpha = 0$, joins them together and puts the resulting vector in decreasing order. More precisely, if $\alpha > 0$, let $B \sim \mathrm{Beta}((1-\alpha)/\alpha, (\theta+\alpha)/\alpha)$, and if $\alpha = 0$, let $B = 1/(\theta+1)$. Let $I_1, I_2, \ldots$ be 0–1 random variables which, given $B$, are independent and identically distributed with $\mathrm{Bernoulli}(B)$ law. Then

$$\mathrm{Coag}_{\alpha,\theta}(x) = \left( \sum_{i \,:\, I_i = 1} x_i, \quad (x_j : I_j = 0) \right)^\downarrow.$$

THEOREM 3.1. *Let $0 \leq \alpha < 1$ and $\theta > -\alpha$. Suppose that $X$ and $Y$ are random variables that take values in $\Delta_\infty^\downarrow$. Then the following statements are equivalent:*

- $X \sim \mathrm{PD}(\alpha, \theta)$ *and, conditional on* $X$, $Y \sim \mathrm{Frag}_\alpha(X)$.
- $Y \sim \mathrm{PD}(\alpha, \theta+1)$ *and, conditional on* $Y$, $X \sim \mathrm{Coag}_{\alpha,\theta}(Y)$.

PROOF. The $\alpha = 0$ case is Proposition 2 of [7] but, for completeness, we reproduce the proof here. Let $(\gamma(t), t \geq 0)$ be a standard gamma process and let $\xi_1 > \xi_2 > \cdots > 0$ be the jumps of $(\gamma(t), 0 \leq t \leq \theta+1)$. Then, by Proposition 2.3,

$$\frac{1}{\gamma(\theta+1)}(\xi_1, \xi_2, \ldots) \sim \mathrm{PD}(0, \theta+1)$$

independently of $\gamma(\theta+1)$, which has a $\mathrm{Gamma}(\theta+1, 1)$ distribution. Now let $\xi_1' > \xi_2' > \cdots > 0$ be the jumps of $(\gamma(t), 0 \leq t \leq 1)$ and let $\xi_1'' > \xi_2'' > \cdots >$ 0 be the jumps of $(\gamma(t), 1 < t \leq \theta + 1)$. Note that coagulating the jumps which happen in the time interval $[0, 1]$ coagulates a proportion $1/(\theta + 1)$ of $\xi_1, \xi_2, \ldots$ chosen uniformly at random. Then the following relationships hold independently:

$$(1) \qquad \frac{1}{\gamma(1)}(\xi_1', \xi_2', \ldots) \sim \mathrm{PD}(0, 1),$$

$$(2) \qquad \frac{1}{\gamma(\theta + 1) - \gamma(1)}(\xi_1'', \xi_2'', \ldots) \sim \mathrm{PD}(0, \theta),$$

$$(3) \qquad \frac{\gamma(1)}{\gamma(\theta + 1)} \sim \mathrm{Beta}(1, \theta).$$

The independence is a consequence of beta–gamma algebra and the fact that the jumps of a subordinator on disjoint time intervals are independent. From (1), it is clear that in the fragmentation step we split with $\mathrm{PD}(0, 1)$. Furthermore, (2), (3) and Proposition 2.2 then imply that

$$\frac{1}{\gamma(\theta + 1)}(\gamma(1), \xi_1'', \xi_2'', \ldots)^{\downarrow} \sim \mathrm{PD}(0, \theta)$$

and that $\gamma(1)/\gamma(\theta + 1)$ is a size-biased pick from this vector.

Suppose now that $\alpha > 0$. Then a variant of the same argument applies. Let $(\tau(t), t \geq 0)$ be a subordinator of Lévy measure $\alpha x^{-\alpha - 1} e^{-x} dx$ and let $S$ be an independent $\mathrm{Gamma}((\theta + 1)/\alpha, \Gamma(1 - \alpha))$ random variable. Let $B$ be an independent $\mathrm{Beta}((1 - \alpha)/\alpha, (\theta + \alpha)/\alpha)$ random variable. Suppose that $\xi_1 > \xi_2 > \cdots > 0$ are the ranked jumps of $(\tau(t), 0 \leq t \leq S)$, that $\xi_1' > \xi_2' > \cdots > 0$ are the ranked jumps of $(\tau(t), 0 \leq t \leq BS)$ and that $\xi_1'' > \xi_2'' > \cdots > 0$ are the ranked jumps of $(\tau(t), BS < t \leq S)$. Then, as $\theta + 1 > 0$, by Proposition 2.4,

$$\frac{1}{\tau(S)}(\xi_1, \xi_2, \ldots) \sim \mathrm{PD}(\alpha, \theta + 1)$$

independently of $\tau(S) = \sum_{i=1}^{\infty} \xi_i$, which has $\mathrm{Gamma}(\theta + 1, 1)$ distribution. Now note that coagulating the jumps that occur in the interval $[0, BS]$ coagulates a proportion $B$ of the jumps $\xi_1, \xi_2, \ldots$. We have $BS \sim \mathrm{Gamma}((1 - \alpha)/\alpha, \Gamma(1 - \alpha))$ by standard beta–gamma algebra. Then the following relationships hold independently:

$$(4) \qquad \frac{1}{\tau(BS)}(\xi_1', \xi_2', \ldots) \sim \mathrm{PD}(\alpha, 1 - \alpha),$$

$$(5) \qquad \frac{1}{\tau(S) - \tau(BS)}(\xi_1'', \xi_2'', \ldots) \sim \mathrm{PD}(\alpha, \theta + \alpha),$$

$$(6) \qquad \frac{\tau(BS)}{\tau(S)} \sim \mathrm{Beta}(1 - \alpha, \theta + \alpha).$$

From (4), we see that in the fragmentation step we split with $\mathrm{PD}(\alpha, 1-\alpha)$. From (5), (6) and Proposition 2.2, we see that

$$\frac{1}{\tau(S)}(\tau(BS), \xi_1'', \xi_2'', \ldots)^{\downarrow} \sim \mathrm{PD}(\alpha, \theta)$$

and that $\frac{\tau(BS)}{\tau(S)}$ is a size-biased pick from this vector. $\square$

REMARKS. (i) Corollary 13 of [18] gives a set of duality relations for certain coagulation and fragmentation operators applied to Poisson–Dirichlet distributions. In particular, for $0 < \alpha < 1$, $0 \le \beta < 1$ and $\theta > -\alpha\beta$, $\mathrm{PD}(\alpha\beta, \theta)$ is fragmented in such a way that each block is split with an independent $\mathrm{PD}(\alpha, -\alpha\beta)$ random variable [call this $(\alpha, -\alpha\beta)$-FRAG]. This results in a $\mathrm{PD}(\alpha, \theta)$ random variable. In reverse, a coagulation of $\mathrm{PD}(\alpha, \theta)$ which coagulates infinitely many different groups of blocks gives $\mathrm{PD}(\alpha\beta, \theta)$ back. The coagulation operator is a little more involved: suppose that the $\mathrm{PD}(\alpha, \theta)$ random variable is $Y = (Y_1, Y_2, \ldots)$. Take an independent random variable, $Q = (Q_1, Q_2, \ldots)$, with $\mathrm{PD}(\beta, \theta/\alpha)$ distribution, and create an open subset $I_Q$ of $[0, 1]$ composed of open intervals whose lengths are given by the vector $Q$:

$$I_Q = (0, Q_1) \cup \bigcup_{i=2}^{\infty} (Q_1 + \cdots + Q_{i-1}, Q_1 + \cdots + Q_i).$$

Now throw independent $\mathrm{U}(0, 1)$ random variables $U_1, U_2, \ldots$ down on the interval. Let $C_1 = \{j \ge 1 : U_j < Q_1\}$ and, for $i \ge 2$, $C_i = \{j \ge 1 : U_j \in (Q_1 + \cdots + Q_{i-1}, Q_1 + \cdots + Q_i)\}$. Finally, let $X = (X_1, X_2, \ldots)$ be obtained by reranking the terms

$$\tilde{X}_i = \sum_{j \in C_i} Y_j, \qquad i \ge 1.$$

Then $X$ has a $\mathrm{PD}(\alpha\beta, \theta)$ distribution. We denote by $(\beta, \theta/\alpha)$-COAG the operation on the vector $Y$ which produces $X$.

Theorem 3.1 provides a companion set of relations (see Figure 1). While Pitman's relations affect the first parameter multiplicatively, our relations affect the second parameter additively.

(ii) For $\alpha = 1/2$ and $\theta = n - 1/2$, $n \ge 1$, Jim Pitman has pointed out that the operation of splitting a size-biased pick from $\mathrm{PD}(\alpha, \theta)$ according to

$$\mathrm{PD}(\alpha\beta, \theta) \quad \xrightarrow{(\alpha, -\alpha\beta)\text{-FRAG}} \quad \mathrm{PD}(\alpha, \theta) \qquad\qquad \mathrm{PD}(\alpha, \theta) \quad \xrightarrow{\text{Frag}_\alpha} \quad \mathrm{PD}(\alpha, \theta+1)$$
$$\xleftarrow[(\beta, \theta/\alpha)\text{-COAG}]{} \qquad\qquad\qquad\qquad \xleftarrow[\text{Coag}_{\alpha,\theta}]{}$$

FIG. 1. Left: *Pitman's duality relations.* Right: *Theorem* 3.1.

PD$(\alpha, 1 - \alpha)$ can be interpreted in terms of the continuum random tree $T$ embedded in a Brownian excursion, as follows. Let $R_n$ be the subtree of $T$ spanned by $n$ points picked at random according to the mass measure $\mu$ of the tree (corresponding to Lebesgue measure on $[0, 1]$). Let $\mu_n$ be the image of $\mu$ on $R_n$ via the map which takes a point $t \in T$ to its closest point in $R_n$. It follows from the line-breaking construction of $R_n$ in [1] that $\mu_n$ is a random discrete distribution whose ranked atoms are distributed according to PD$(1/2, n - 1/2)$ and that, in the growth step from $R_n$ to $R_{n+1}$, a size-biased choice of one of these atoms is split according to PD$(1/2, 1/2)$ to create the atoms of $\mu_{n+1}$. The inverse coagulation operation can also be seen in this setting as a corollary of Aldous's results. It appears that similar interpretations in terms of continuum trees for other values of $(\alpha, \theta)$ can be based on Section 5.3 of [12]. Such interpretations are the subject of work in progress by Pitman and Winkel. (See [2, 3, 19] for further background.)

## 4. Fragmentation and coagulation processes.

Define a discrete-time Markov fragmentation chain $(X(i), i \geq 0)$ that takes values in $\Delta_\infty^{\downarrow}$ as follows: for $i \geq 0$, conditional on $X(i)$,

$$X(i + 1) \sim \mathrm{Frag}_\alpha(X(i)).$$

If $X(0) \sim$ PD$(\alpha, \theta)$, then, by Theorem 3.1, $X(i) \sim$ PD$(\alpha, \theta + i)$ for $i \geq 1$.

Likewise, define an inhomogeneous Markov coagulation chain by

$$\ldots, X(i + 1), X(i), \ldots, X(1), X(0).$$

By Theorem 3.1, conditional on $X(i + 1)$,

$$X(i) \sim \mathrm{Coag}_{\alpha, \theta + i}(X(i + 1))$$

for $i \geq 0$.

We can construct a continuous-time Markov fragmentation process $(Y(t), t \geq 0)$ by taking an independent standard Poisson process $(N(t), t \geq 0)$ and letting

$$Y(t) = X(N(t))$$

for $t \geq 0$. In the terminology of Bertoin [6], this is a self-similar fragmentation of index 1, erosion coefficient 0 and dislocation measure PD$(\alpha, 1 - \alpha)$. Since the dislocation measure is finite (it is a probability measure), the fact that the fragmentation has index of self-similarity $\delta$ just means that a block of size $x$ splits at a rate proportional to $x^\delta$, and independently of the other blocks. Since here the total rate of splitting is 1 at any time, and we split a size-biased pick from among the blocks, we must have each block splitting at the rate of its length, that is, $\delta = 1$.

Suppose that we fix $\theta$. In the $\alpha = 0$ case, it was shown in [7] that a dual Markovian coagulation chain may be defined by $(Y(e^{-t}), t \geq 0)$. It can be checked that when the state has distribution $\mathrm{PD}(0, \theta + n)$, for any $n \geq 1$, the process waits an exponential time of parameter $n$ and then jumps to a state distributed as $\mathrm{PD}(0, \theta + n - 1)$. In principle, the same construction may be performed in the $\alpha > 0$ case. However, here the inhomogeneity of the discrete-time coagulation chain becomes a problem. In the $\alpha = 0$ case, the distributions $\mathrm{PD}(0, \theta)$ are almost surely distinguishable as $\theta$ varies. Thus, it is possible to tell from the current state which coagulation operator to apply to it to get the next state. In the $\alpha > 0$ case, however, the distributions $\mathrm{PD}(\alpha, \theta)$ are mutually absolutely continuous as $\theta$ varies. Thus, it is not possible to detect almost surely from the state what the second parameter is and then work out which coagulation operator to apply.

REMARKS. In [7], the $\alpha = 0$ case of these processes is shown to arise naturally in the context of the genealogy of certain continuous-state branching processes. We do not see any way to generalize those results to the case $\alpha > 0$.

## 5. Random recursive trees.
Let $\alpha \in [0, 1)$ and $\theta > -\alpha$. An *exchangeable* $(\alpha, \theta)$ *partition* of $\mathbb{N}$ (or of any infinite subset $A \subseteq \mathbb{N}$) is defined as follows:

(i) Generate a $\mathrm{PD}(\alpha, \theta)$-distributed vector $(Y_1, Y_2, Y_3, \ldots)$.

(ii) Conditionally on $Y_1, Y_2, Y_3, \ldots$, assign $i$ to block $j$ with probability $Y_j$, independently for each $i \in A$.

Suppose we order the blocks of this partition in increasing order of their smallest elements. Let $B_i$ be the $i$th block and let $F_i$ be its *asymptotic frequency*, that is,

$$F_i = \lim_{n \to \infty} \frac{|B_i \cap \{1, 2, \ldots, n\}|}{|A \cap \{1, 2, \ldots, n\}|}.$$

Then $(F_1, F_2, \ldots)$ is a size-biased ordering of $(Y_1, Y_2, \ldots)$. In particular, $(F_1, F_2, \ldots)$ has the $\mathrm{GEM}(\alpha, \theta)$ distribution.

An alternative way to construct an exchangeable $(\alpha, \theta)$ partition is via the *Chinese restaurant process* of Dubins and Pitman (see [19]), defined as follows. Person 1 enters a Chinese restaurant and sits at the first table. Person 2 sits either at the same table or at a new one; in general, each subsequent person (numbered successively $3, 4, \ldots$) sits either at one of the occupied tables or at a new one. Suppose that people $1, 2, \ldots, n$ have sat at $k$ tables, where table $i$ has $n_i$ customers (with $n_i \geq 1$ and $\sum_{i=1}^{k} n_i = n$). Then person $n + 1$ starts a new table with probability

$$\frac{\theta + k\alpha}{n + \theta}$$

and sits at table $i$ with probability

$$\frac{n_i - \alpha}{n + \theta}$$

for $1 \leq i \leq k$ (see Figure 2). The partition of $\mathbb{N}$ into blocks that correspond to the different tables is an exchangeable $(\alpha, \theta)$ partition of $\mathbb{N}$. (Of course, the construction for general $A \subseteq \mathbb{N}$ is analogous.)

We are now ready to describe the model of an $(\alpha, \theta)$-*recursive tree*. A recursive tree is a rooted labeled tree such that the vertex labels increase along paths away from the root. We now construct a random recursive tree as follows. We start with the root labeled 0 and a single child labeled 1. Vertices $2, 3, \ldots$ are now added in turn; vertex $i$ is added as a child of one of the existing vertices $0, 1, \ldots, i - 1$. When vertex $i$ is added, the probability that it is added as a child of vertex $j$ is proportional to $1 - \alpha + \alpha k_j$ if $j \geq 1$, where $k_j$ is the current number of children of vertex $j$. It is added as a child of vertex 0 with probability proportional to $\theta + \alpha k_0$. A tree on $\mathbb{Z}_+$ that arises in this way is called an $(\alpha, \theta)$-recursive tree.

It is useful to generate the same distribution via a continuous-time Markov chain. Namely, we again start with the root labeled 0 and a single child labeled 1. From then on, new children of the root arrive at rate $\theta + \alpha k_0$ and new children of vertex $j$, $j \geq 1$, arrive at rate $1 - \alpha + \alpha k_j$ (where, again, $k_j$ is the current number of children of vertex $j$). The new vertices are numbered in the order they arrive. This continuous-time construction will make it possible to deduce directly various useful independence properties of the tree.

Now consider the following procedure. Remove vertex 0 and record the partition of $\mathbb{N}$ given by the resulting forest. Call this partition $B^{(0)}$, where the blocks $B_1^{(0)}, B_2^{(0)}, \ldots$ are listed in increasing order of smallest element (note that this smallest element is necessarily the root of one of the recursive trees in the forest). Now, for $i \geq 1$, define $B^{(i)}$ to be the partition of $\mathbb{N} \setminus \{1, 2, \ldots, i\}$ obtained by removing vertices $0, 1, 2, \ldots, i$ (again, the blocks $B_1^{(i)}, B_2^{(i)}, \ldots$ of this partition are listed in increasing order of smallest element).

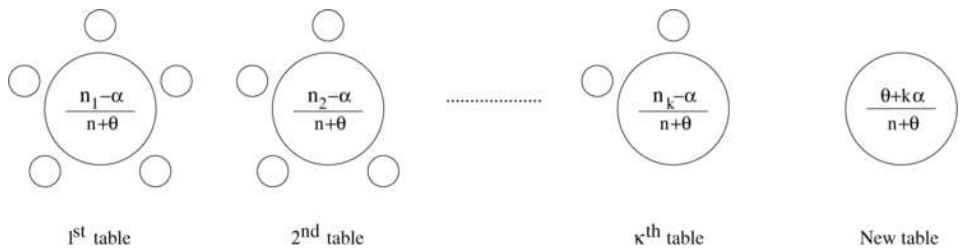

Fig. 2. *The Chinese restaurant process.*

THEOREM 5.1. *For all $i \geq 0$, the blocks $B_1^{(i)}, B_2^{(i)}, \ldots$ form an exchangeable $(\alpha, \theta + i)$ partition of $\mathbb{N} \setminus \{1, 2, \ldots, i\}$. In particular, $B_1^{(i)}, B_2^{(i)}, \ldots$ possess asymptotic frequencies $F_1^{(i)}, F_2^{(i)}, \ldots$ such that $(F^{(i)})^{\downarrow} \in \Delta_{\infty}^{\downarrow}$ and $F^{(i)} \sim \mathrm{GEM}(\alpha, \theta + i)$ for all $i \geq 0$. Moreover, letting $G^{(i)} = (F^{(i)})^{\downarrow}$, we have that $(G^{(i)}, i \geq 0)$ is a Markov chain such that the following statements hold:*

- *For $i \geq 1$, conditional on $G^{(i)}$, $G^{(i+1)}$ has the same distribution as $\mathrm{Frag}_{\alpha}(G^{(i)})$.*
- *For $i \geq 1$, conditional on $G^{(i+1)}$, $G^{(i)}$ has the same distribution as $\mathrm{Coag}_{\alpha, \theta+i}(G^{(i+1)})$.*

This result is illustrated in Figure 3. It entails that as long as $X(0) \sim \mathrm{PD}(\alpha, \theta)$, then

$$(G^{(i)}, i \geq 0) \overset{d}{=} (X(i), i \geq 0),$$

where $(X(i), i \geq 0)$ is the Markov fragmentation chain of Section 4, and clearly implies Theorem 3.1. The size-biased view of the fragmentation chain given by the random recursive tree seems a very natural description: rather than having two sources of external randomness (one to take a size-biased pick from the state vector and another to split it), here the randomness is entirely in the tree; given the tree, the fragmentation is deterministic. The tree can be thought of as a concrete representation of the filtration of the fragmentation.

REMARKS. (i) For a general survey of results on recursive trees, see [22]. Random recursive trees similar to certain $(\alpha, \theta)$-recursive trees have been studied by Szymański [23] and by Mahmoud, Smythe and Szymański [14]. One object of interest is the *k*th *branch*, that is, the subtree rooted at $k$ in a random recursive tree labeled by $\{0, 1, \ldots, n\}$. Call the size of the *k*th branch $T_{n,k}$. In our model, $F_1^{(k-1)}$ is the almost sure limit of $T_{n,k}/n$ as $n \to \infty$. In his Theorem 8, Szymański [23] finds the mean and variance of $T_{n,k}$ in the $(1/2, 0)$-recursive tree; his results are consistent with the limiting mean and variance implied by Theorem 5.1. Theorem 5 of [14] gives the limiting distribution of $T_{n,k}/n$ as $\mathrm{Beta}(1/2, k)$ in a model which is essentially the $(1/2, 1/2)$-recursive tree; this is also what we expect from Theorem 5.1. Various related classes of random graphs and trees constructed by preferential attachment are considered, for example, by Barabási and Albert [5], Bollobás, Riordan, Spencer and Tusnády [9], Móri [15] and Rudas, Tóth and Valkó [21].

(ii) Taking $\alpha = -1/m$ and $\theta = r/m$ for integers $m \geq 1$ and $r \geq 2$ in the Chinese restaurant process gives an exchangeable random partition into $r$ blocks whose asymptotic frequencies are the decreasing rearrangement of a

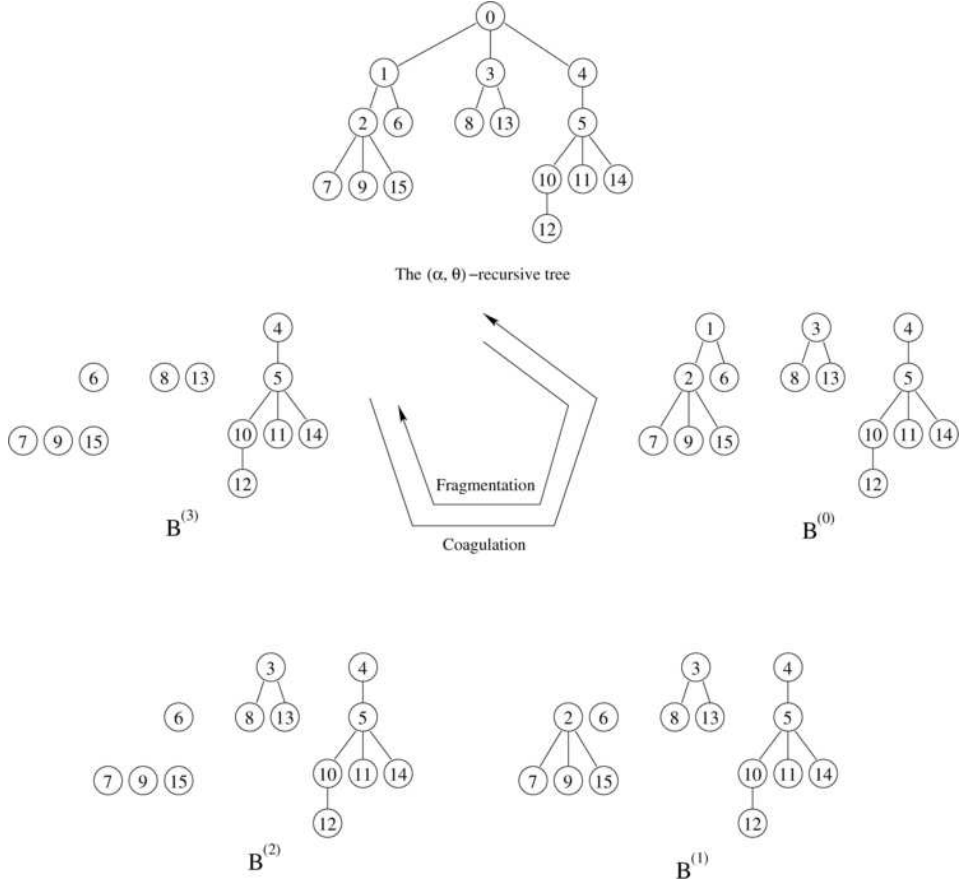

The $(\alpha, \theta)$–recursive tree

$\mathbf{B}^{(3)}$

Fragmentation

Coagulation

$\mathbf{B}^{(0)}$

$\mathbf{B}^{(2)}$

$\mathbf{B}^{(1)}$

Fig. 3. *Fragmentation and coagulation for the $(\alpha, \theta)$-recursive tree.*

Dirichlet$(1/m, \ldots, 1/m)$ random vector with $r$ parameters all equal to $1/m$. The construction of the $(\alpha, \theta)$-recursive tree works for these parameters as well [to be precise, the probability we add a vertex to the root is proportional to $(r - k_0)_+/m$ and the probability we add to any other vertex is proportional to $(m + 1 - k_j)_+/m$], giving a tree whose out-degrees are all equal to $m + 1$, except for the root which has degree $r$. Analogs of Theorems 3.1 and 5.1 hold in this case too; see [7] for details.

**6. Proof of Theorem 5.1.** We prove Theorem 5.1 in a series of lemmas. (Note that this proof is independent of Theorem 3.1.)

LEMMA 6.1. *The blocks $B_1^{(i)}, B_2^{(i)}, \ldots$ form an exchangeable $(\alpha, \theta + i)$ partition of $\mathbb{N} \setminus \{1, \ldots, i\}$.*

Proof.    By the children of a *set* of vertices, we will mean here vertices which are neighbors of that set in the tree, but are not contained within it. We work via the Chinese restaurant process. Think of children of the set of vertices $\{0, 1, 2, \ldots, i\}$ as starting new tables. Consider the construction of the $(\alpha, \theta)$-recursive tree and suppose that vertices labeled $i+1, i+2, \ldots, i+n$ have already arrived, including $k$ children of the set $\{0, 1, 2, \ldots, i\}$. Suppose that the $k$ subtrees rooted at these children have sizes $n_1, n_2, \ldots, n_k$, where $\sum_{j=1}^{k} n_j = n$ (the subtrees correspond to tables in the Chinese restaurant process). Then $i + n + 1$ forms a new table at rate $\theta + i(1 - \alpha) + \alpha(k + i) = \theta + i + \alpha k$ and is added to table $j$ at rate $n_j(1 - \alpha) + \alpha(n_j - 1) = n_j - \alpha$ for $1 \leq j \leq k$. Hence, the total rate is $\sum_{j=1}^{k} n_j - k\alpha + \theta + i + \alpha k = n + \theta + i$. Thus, $i + n + 1$ forms a new table with probability $(\theta + i + \alpha k)/(n + \theta + i)$ and adds to table $j$ with probability $(n_j - \alpha)/(n + \theta + i)$ for $1 \leq j \leq k$. Hence, we have a Chinese restaurant process of parameters $(\alpha, \theta + i)$ and so removing vertices $0, 1, 2, \ldots, i$ gives blocks which form an $(\alpha, \theta + i)$ partition.    □

We now describe the $(\alpha, \theta)$-recursive tree in terms of a set of "nested" $(\alpha, \theta)$ and $(\alpha, 1 - \alpha)$ partitions.

Lemma 6.2.    *The distribution of the $(\alpha, \theta)$-recursive tree is characterized (among distributions on recursive trees on $\mathbb{Z}_+$) by the following properties:*

(i) *The blocks that correspond to subtrees rooted at different children of the root form an exchangeable $(\alpha, \theta)$ partition of $\mathbb{N}$.*

(ii) *Let $i \geq 1$. Consider the set $D(i) \subseteq \mathbb{N} \setminus \{1, 2, \ldots, i\}$ of labels which are descendants of $i$ in the tree. Conditional on $D(i)$, the blocks that correspond to subtrees rooted at different children of $i$ form an exchangeable $(\alpha, 1 - \alpha)$ partition of $D(i)$ and this partition is independent of the structure of the subtree on $\mathbb{Z}_+ \setminus D(i)$.*

One could say that the subtree rooted at any vertex except the root has the structure of an $(\alpha, 1 - \alpha)$-recursive tree; in the special case $\theta = 1 - \alpha$, the whole tree also has this structure and so one has full self-similarity.

Proof of Lemma 6.2.    The partition of $\mathbb{N}$ into blocks as in part (i) determines the labels of the children of the root (because of the property that labels increase along paths away from the root, so that the children of the root are the smallest elements of each block of the partition). Then the subpartitions of the blocks (minus their smallest elements) as in part (ii) determine the next level of the tree and so on; thus we indeed have a characterization of the distribution of the tree.

Property (i) is the $i = 0$ case of Lemma 6.1. The self-similarity and independence properties in (ii) are most easily seen via the continuous-time construction of the tree. Considering only the subtree rooted at vertex $i$, we see that the process of vertices arriving at that subtree corresponds precisely to that for building an $(\alpha, 1-\alpha)$ tree (with a different set of labels); its structure depends on the evolution of the rest of the tree only through the sequence of labels which are assigned to new vertices that join the subtree. □

The following alternative procedure also builds the $(\alpha, \theta)$-recursive tree and in addition will describe precisely the evolution of the process $(G^{(i)}, i \geq 0)$ obtained by removing the vertices $0, 1, 2, \ldots$ in turn. This construction employs the self-similarity and naturally gives the Markov property of the fragmentation process, which is not easily obtained by other means.

At stage $n$ of the procedure we have a tree with $n + 1$ internal vertices labeled $0, 1, 2, \ldots, n$, each of which has infinitely many children which are leaves. At each stage we label one of the leaves and create new leaves which are children of the newly labeled vertex (which, of course, ceases to be a leaf itself). Each vertex of the tree carries a weight. The weight of a vertex represents the asymptotic frequency of the set of all the labels which will be assigned to the descendants of that vertex. In particular:

- The weight of the root is 1.
- For any internal (i.e., already labeled) vertex, the weight of the vertex equals the sum of the weights of all its children. [In fact, the weights of the children are a splitting, by a $\mathrm{PD}(\alpha, 1-\alpha)$ random vector, of the weight of the parent, unless the parent is the root, in which case the splitting is by $\mathrm{PD}(\alpha, \theta)$.]

We start at stage 0 with a tree that consists of the root, labeled 0, with weight 1, and an infinite number of unlabeled leaves which are children of the root and whose weights are given by a $\mathrm{PD}(\alpha, \theta)$ random vector.

To pass to stage 1, we choose one of the children of the root in a size-biased way (according to the weights assigned to the vertices) and assign label 1 to this vertex. We create an infinite number of children of the newly labeled vertex (which become leaves); these are assigned weights given by a splitting of the weight of the parent vertex by a $\mathrm{PD}(\alpha, 1-\alpha)$ random vector.

To pass from stage $n$ to stage $n+1$, we similarly choose between the children of the root in a size-biased way. If we choose an already labeled vertex, we now choose between the children of that chosen vertex in a size-biased way. This continues until we reach a leaf. This leaf is now labeled $n+1$; its children are created and their weights are assigned by a $\mathrm{PD}(\alpha, 1-\alpha)$ splitting as above.

All of the size-biased picks and random vectors are independent.

Proceeding in this way, one constructs a tree on $0, 1, 2, \ldots$. From the description of the random recursive tree in Lemma 6.2 in terms of nested exchangeable $(\alpha, \theta)$ and $(\alpha, 1 - \alpha)$ partitions, it follows that this tree indeed has the distribution of an $(\alpha, \theta)$-recursive tree.

We can now use this construction to identify the chain $(G^{(i)}, i \geq 0)$ with the chain obtained by applying $\text{Frag}_\alpha$ repeatedly.

LEMMA 6.3.  *The process* $(G^{(i)}, i \geq 0)$ *is a Markov chain such that conditional on* $G^{(i)}$, $G^{(i+1)}$ *has the same distribution as* $\text{Frag}_\alpha(G^{(i)})$.

PROOF.  First note that the weights of all the leaves at stage $i$ correspond to the state $G^{(i)}$, recording the asymptotic frequencies of the subtrees obtained when the vertices $0, 1, \ldots, i$ are removed and ranking these frequencies in decreasing order. The procedure of passing from stage $i$ to $i + 1$—choosing first between the children of the root in a size-biased way, then between the children of that child in a size-biased way and so on until a leaf is chosen—is equivalent simply to choosing between all the leaves in a size-biased way. Given the state of the tree at stage $i$, this choice of vertex $i + 1$ is made independently of previous choices. Thus, given $G^{(i)}$, we obtain $G^{(i+1)}$ precisely by applying the random operator $\text{Frag}_\alpha$ to $G^{(i)}$ independently of $G^{(0)}, \ldots, G^{(i-1)}$.  □

We have proved that $(G^{(i)}, i \geq 0)$ is a Markov chain and that $G^{(i+1)}$ is the required fragmentation of $G^{(i)}$. It remains only to show the coagulation property. To understand the coagulation mechanism in the tree, we need some more notation. For $i \geq 0$, starting from blocks $B_1^{(i+1)}, B_2^{(i+1)}, \ldots$, create new blocks $\tilde{B}_1^{(i)}, \tilde{B}_2^{(i)}, \ldots$ as follows. Let $\tilde{I}_1^{(i+1)}$ be an independent Bernoulli$((1 - \alpha)/(\theta + i + 1))$ random variable and recursively generate $\tilde{I}_2^{(i+1)}, \tilde{I}_3^{(i+1)}, \ldots$ that take values 0 or 1 via

$$(7) \qquad \mathbb{P}(\tilde{I}_{k+1}^{(i+1)} = 1 | \tilde{I}_1^{(i+1)}, \tilde{I}_2^{(i+1)}, \ldots, \tilde{I}_k^{(i+1)}) = \frac{1 - \alpha + \alpha \sum_{j=1}^k \tilde{I}_j^{(i+1)}}{\theta + i + 1 + \alpha k}$$

independently of $B^{(i+1)}$. Let $C^{(i+1)} = \{k : \tilde{I}_k^{(i+1)} = 1\}$, let

$$\tilde{B}_1^{(i)} = \{i + 1\} \cup \bigcup_{k \in C^{(i+1)}} B_k^{(i+1)}$$

and let $\tilde{B}_2^{(i)}, \tilde{B}_3^{(i)}, \ldots$ be the blocks $B_k^{(i+1)}$ not contained in the union $\tilde{B}_1^{(i)}$, listed in increasing order of smallest element (notice that this retains the ordering of the indices from the previous step). Let $\tilde{B}^{(i)} = (\tilde{B}_1^{(i)}, \tilde{B}_2^{(i)}, \ldots)$.

LEMMA 6.4.  *For* $i \geq 1$, $\tilde{B}^{(i)}$ *is a partition of* $\mathbb{N} \setminus \{1, 2, \ldots, i\}$ *and*

$$(\tilde{B}^{(i)}, B^{(i+1)}) \stackrel{d}{=} (B^{(i)}, B^{(i+1)}).$$

PROOF. Suppose that we are given $B^{(i+1)}$. We know that in the full recursive tree, each of the blocks $B_1^{(i+1)}, B_2^{(i+1)}, \ldots$ is connected to one of the vertices $0, 1, \ldots, i+1$. We need to know which of them are connected to $i+1$.

Let $(a_1, a_2, a_3, \ldots)$ be the sequence of children of the nodes $\{0, 1, \ldots, i+1\}$ in order. (In particular, $a_1 = i+2$.) Let $(b_1, b_2, \ldots)$ be $\{i+2, i+3, \ldots\} \setminus \{a_1, a_2 \ldots\}$ in order. Now for $j \geq 1$, let $p_j$ be the parent of node $a_j$ and let $q_j$ be the parent of node $b_j$. (Thus $p_j \in \{0, 1, \ldots, i+1\}$ and $q_j \in \{i+2, i+3, \ldots\}$.) From the continuous-time construction of the random recursive tree, one can deduce that the three sequences $(a_1, a_2, \ldots)$, $(p_1, p_2, \ldots)$ and $(q_1, q_2, \ldots)$ are independent. Since $B^{(i+1)}$ is a function of $(a_1, a_2, \ldots)$ and $(q_1, q_2, \ldots)$, we therefore have that $(p_1, p_2, \ldots)$ is independent of $B^{(i+1)}$.

Since $p_j$ is the vertex to which block $B_j^{(i+1)}$ is attached, we need to know specifically which of the $p_j$ are equal to $i+1$. Now, regardless of the form of the subtree spanned by vertices $0, 1, \ldots, i+1$, vertex $i+1$ has children at rate $1 - \alpha + \alpha n_{i+1}$, where $n_{i+1}$ is the number of children that it has already, and the group $0, 1, \ldots, i$ has children at total rate $\theta + i + \alpha + \alpha \tilde{n}$, where $\tilde{n}$ is the combined total number of children that the group already has other than $i+1$. Hence, the probability that $p_1 = i+1$ is equal to $(1-\alpha)/(\theta+i+1)$. For $k \geq 1$, let $I_k^{(i+1)}$ be the indicator function of the event $\{p_k = i+1\}$, which of course is the same as the event $\{B_k^{(i+1)} \text{ is attached to } i+1\}$. Then for $k \geq 1$ we have

$$\mathbb{P}(I_{k+1}^{(i+1)} = 1 | I_1^{(i+1)}, I_2^{(i+1)}, \ldots, I_k^{(i+1)}) = \frac{1 - \alpha + \alpha \sum_{j=1}^k I_j^{(i+1)}}{\theta + i + 1 + \alpha k},$$

so that the sequence $(I_1^{(i+1)}, I_2^{(i+1)}, \ldots)$ has the same law as the sequence $(\tilde{I}_1^{(i+1)}, \tilde{I}_2^{(i+1)}, \ldots)$ constructed at (7). Thus, conditional on $B^{(i+1)}$, $\tilde{B}^{(i)} \stackrel{d}{=} B^{(i)}$. $\square$

The proof of Theorem 5.1 is now completed by the following lemma.

LEMMA 6.5. *Conditional on $G^{(i+1)}$, $G^{(i)}$ has the same distribution as* $\mathrm{Coag}_{\alpha, \theta+i}(G^{(i+1)})$.

PROOF. Consider the construction of $\tilde{B}^{(i)}$ from $B^{(i+1)}$. The random variables $\tilde{I}_1^{(i+1)}, \tilde{I}_2^{(i+1)}, \ldots$ defined at (7) describe the evolution of a "generalized Pólya urn." In particular, the limit

$$B = \lim_{k \to \infty} \frac{1}{k} \sum_{j=1}^k \tilde{I}_j^{(i+1)}$$

exists almost surely and has a Beta$((1 - \alpha)/\alpha, (\theta + i + \alpha)/\alpha)$ distribution [except in the case $\alpha = 0$, when $B = 1/(\theta + i + 1)$ a.s.]; then, conditional on $B$, the variables $\tilde{I}_1^{(i+1)}, \tilde{I}_2^{(i+1)}, \ldots$ are independent and identically distributed Bernoulli$(B)$ random variables (see [8] or [17]).

Thus, conditional on $B^{(i+1)}$, the ranked asymptotic frequencies $\tilde{G}^{(i)}$ of $\tilde{B}^{(i)}$ have the same distribution as $\mathrm{Coag}_{\alpha, \theta+i}(G^{(i+1)})$. Since this distribution depends on $B^{(i+1)}$ only through $G^{(i+1)}$, we have that the same is true if we condition instead on $G^{(i+1)}$ (this follows by Dynkin's criterion for a function of a Markov process to be Markov; see Theorem 10.13 of [11]). But by Lemma 6.4, we have $(\tilde{B}^{(i)}, B^{(i+1)}) \stackrel{d}{=} (B^{(i)}, B^{(i+1)})$, so also $(\tilde{G}^{(i)}, G^{(i+1)}) \stackrel{d}{=} (G^{(i)}, G^{(i+1)})$. Hence, conditional on $G^{(i+1)}$, $G^{(i)}$ has the same distribution as $\mathrm{Coag}_{\alpha, \theta+i}(G^{(i+1)})$.   □

## 7. Continuous-time branching models.

Pitman [16], gave a construction of the two-parameter partition structure [or, equivalently, the $(\alpha, \theta)$ Chinese restaurant process] via continuous-time branching models (possibly with immigration). In this section we discuss the relationship between such a construction and our recursive tree model.

We first describe Pitman's construction. Let $0 \leq \alpha < 1$, $\theta > -\alpha$. Consider a population of individuals of two types: novel and clone. Each individual is assigned a color: a novel individual always has a new color, while a clone is always assigned the same color as its parent. Every individual has an infinite lifetime. Starting from a single novel individual at time $t = 0$, this first individual produces novel offspring according to a Poisson process of rate $\theta + \alpha$ and clone offspring according to an independent Poisson process of rate $1 - \alpha$. The reproduction rules for the other individuals are as follows:

- Novel individuals produce novel offspring according to a Poisson process of rate $\alpha$ and independently produce clone offspring according to a Poisson process with rate $1 - \alpha$.
- Clone individuals produce clone offspring according to a Poisson process of rate 1.

Individuals are labeled in the order they appear. The colors of individuals naturally induce a random partition of $\mathbb{N}$, which has the PD$(\alpha, \theta)$ distribution, by comparison of the growth procedure with the $(\alpha, \theta)$ Chinese restaurant process.

If $\theta > 0$, we may treat the first individual just like any other novel individual (i.e., producing novel individuals at rate $\alpha$ and clones at rate $1 - \alpha$) and introduce an independent Poisson migration process of novel individuals which arrive at rate $\theta$. This way of looking at things provides an easy way to see the fact (due to Pitman [18]) that taking a PD$(0, \theta)$ random variable and splitting each block with an independent PD$(\alpha, 0)$ random variable gives a PD($\alpha, \theta$) random variable. Indeed, by ignoring the clone/novel difference and just looking at the partition generated by the descendencies of the immigrant individuals, we see a $(0, \theta)$ partition. If we then keep track of the different colors as well, we see a refinement of this partition. Moreover, within each block of the coarse partition, the colors are generated according to the rules for an $(\alpha, 0)$ partition.

Let us explain where the branching model construction of an $(\alpha, \theta)$ partition (in the no immigration setting) differs from the recursive tree construction. Novel individuals in the branching model are exactly the children of the root in the recursive tree (although which novel individual is the "child" of which other novel individual has no meaning in the recursive tree setup). The growth rates for the subfamilies are the same: the collection consisting of a novel individual and its clone descendents (say $k$ of them) produces new clone individuals at rate $1 - \alpha + k$. In the recursive tree, likewise the subtree descending from a particular child of the root grows at rate $(1 - \alpha)(k + 1) + \alpha k = 1 - \alpha + k$. However, the new individuals are not being added in the same places on the subtree. For example, a novel individual in Pitman's construction has clone children at fixed rate $1 - \alpha$. In our subtree, however, the individual at the top of the subtree has children at rate $1 - \alpha + \alpha \# \{$children it has already had$\}$. These genealogical differences make no difference to the partition obtained, but do make a difference to the nesting property for successive partitions.

The coagulation–fragmentation duality (Theorem 3.1) for $0 \leq \alpha < 1$, $\theta > -\alpha$ can be interpreted through a variant of the branching model construction. We introduce killing, as in the recursive tree, and also give the clones hidden features. These hidden features enable us to obtain the required nesting of partitions. More precisely:

- Each clone individual is different from its brothers and has a brand new color, but this difference (in type as well as in color) is invisible until its parent is killed.
- Each clone individual actually generates novel individuals at rate $\alpha$ and independently generates clone individuals at rate $1 - \alpha$, but this difference (in type as well as in color) among its offspring is invisible until its parent is killed and it becomes novel.

As in the recursive trees setting, start with a PD($\alpha, \theta$) population and kill the first individual, then the second, then the third and so forth. Notice that immediately after killing an individual, its clone children and the offspring of those clone children may change their novel/clone status or color and form a number of new blocks which fragment the block of the killed individual. This provides an alternative encoding of the fragmentation chain.

**Acknowledgments.**    Part of this work was carried out during the program on Probability, Algorithms and Statistical Physics at the MSRI in Berkeley in spring 2005. Christina Goldschmidt and James Martin would like to thank the MSRI for its support. We are very grateful to Jim Pitman for valuable discussions.

R. Dong
Department of Statistics
367 Evans Hall
University of California
Berkeley, California 94720-3860
USA
E-mail: ruidong@stat.berkeley.edu

C. Goldschmidt
Statistical Laboratory
Centre for Mathematical Sciences
Cambridge University
Wilberforce Road
Cambridge CB3 0WB
United Kingdom
E-mail: C.Goldschmidt@statslab.cam.ac.uk

J. B. Martin
Department of Statistics
Oxford University
1 South Parks Road
Oxford OX1 3TG
United Kingdom
E-mail: martin@stats.ox.ac.uk